\documentclass[12pt]{article}
\usepackage{amsmath,amsfonts,amssymb,layout}
\begin{document}
\def\B{{\cal B}}
\def\O{{\cal O}}
\def\P{{\Bbb P}}
\def\Q{{\Bbb Q}}
\def\C{{\Bbb C}}
\def\Z{{\Bbb Z}}
\def\qed{\hfill\vbox{\hrule\hbox{\vrule\kern3pt\vbox{\kern6pt}\kern3pt\vrule}\hrule}\bigskip}
\newcommand\Bt{{\B_{\mathrm t}}}
\newcommand\Be{{\B_{\mathrm e}}}
\newcommand\Proj{\mathrm{Proj}}
\newcommand\Spec{\mathrm{Spec}}

\def\qed{\hfill\vbox{\hrule\hbox{\vrule\kern3pt\vbox{\kern6pt}\kern3pt\vrule}\hrule}\bigskip}

%\baselineskip = 16pt
%\font\FN=cmr10 at 18pt

%\theoremstyle{definition}
\newtheorem{prop}{\bf Proposition}[section]
\newtheorem{tma}[prop]{Theorem}
\newtheorem{lma}[prop]{Lemma}
\newtheorem{cor}[prop]{Corollary}
\newtheorem{MAIN}[prop]{THEOREM}
\newtheorem{rem}[prop]{Remark}
\newtheorem{ex}[prop]{Example}
\newenvironment{pf}{{\it Proof.}}{\qed}
\newenvironment{df}{\bigskip\noindent{\bf Definition}}{\bigskip}
\newenvironment{TMA}{\bigskip\noindent{THEOREM}\it}{\bigskip}

%\section{Noncyclic torsion}

\title{Fano threefolds with noncyclic torsion in the divisor class group}

\author{Jorge Caravantes}

\maketitle

\begin{abstract}
In this note we study Fano threefolds with noncyclic torsion in the divisor class group. Since they can all be obtained as quotients of Fano threefolds, we get also all examples that can be obtained as quotients of low codimension Fanos in the weighted projective space.
\end{abstract}

Reid's graded rings method has been used to find families of examples of manifolds in the weighted projective space. If we restrict to Fano threefolds,  there are lists of codimension 1, 2 and 3 Fano threefolds in a weighted projective space, due to Reid, Fletcher and Alt{\i}nok respectively, and a forthcoming list in codimension 4 is almost finished. All examples of these and other lists can be found in \cite{GRDW}. Recently, the author has written in \cite{C} all quotients (with at most cyclic quotent terminal singularities) by a finite cyclic group of all codimension 1, 2 and 3 Fano threefolds in a weighted projective space.

In this note we continue the work in \cite{C} and take the problem of the existence of noncyclic torsion in the divisor class group of a $\Q$-Fano variety $X$. Our main result is the following:

\begin{TMA}
All possible Fano-Enriques quotients of a codimension 1, 2 or 3 Fano threefold by a finite noncyclic group action are those listed in Section \ref{listas}.
\end{TMA}

The proof is a case by case one, so we will just see first the method and then develop it in detail in only one example in this paper since we do not think it is useful to write down all cases. In the first section we restrict to all possible pairs of orders of independent torsion divisors that can be found in a Fano threefold. In Section 2 we find all possible torsion baskets for pairs of torsion divisors and study the cases of three and four independent torsion divisors. In Section 3 we develop the method to find noncyclic Fano--Enriques quotients from Fano threefolds. Finally we list all low codimension examples in Section 4.

\section{First results on the torsion group}

In this section we determine all possible forms that the torsion subgroup of the divisor class group of a Fano threefold $X$ can take. Let us suppose that $\sigma$ and $\tau$ are two independent torsion divisors on our Fano--Enriques threefold $X$ of orders $r$ and $s$ respectively. By the structure theorem of finitely generated abelian groups, we can suppose that $r$ divides $s$. We admit for our varieties only terminal cyclic quotient singularites. This means that any singularity ``Q" in this paper is, locally, the quotient of $\C^3$ by an action $(x,y,z)\mapsto (e^{\frac{2\pi i}{r_Q}}x,e^{a\frac{2\pi i}{r_Q}}y,e^{-a\frac{2\pi i}{r_Q}}z)$ where $a$ and $r_Q$ are coprime. Now, if $\sigma$ and $\tau$ are $l_Q{K_X}$ and $m_Q{K_X}$ respectively in a sufficiently small analytic neighbourhood of $Q$ (so $0\le l_Q,m_Q\le r_Q-1$), we denote the singularity type this way:
$$ \frac{1}{r_Q}(1,a,-a)_{l_Q,m_Q},$$
suppressing any vanishing subindex. We recall from \cite{C} that $l_Q,m_Q\in \{2,3,4,5,6,8\}$.

Now, we can define
$$\Bt(\sigma):=\{ Q\ |\ l_Q\ne 0 \}$$
$$\Bt(\tau):=\{ Q\ |\ m_Q\ne 0 \}$$
$$\Bt:=\Bt(\sigma)\cup\Bt(\tau)$$
$$\Be:=\{ Q\in \mathrm{Sing} X\ |\ l_Q=m_Q=0\}$$
and we have the following properties:

\begin{rem}
We can cosider any linear combination of $\sigma$ and $\tau$ separately, so we can recall from \cite{C} that
\begin{enumerate}
\item $\Bt(\sigma)$ and $\Bt(\tau)$ are among the torsion subsets in Table 1 in \cite{C}, each one appearing in the box refering to the appropriate torsion ($r$ and $s$ respectively).
\item if $\Bt(\sigma+\tau)$ is, coherently, $\{ Q | l_Q+m_Q \not \sim 0\ \mbox{mod}\ r_Q \}$ then $\B(\sigma+\tau)$ is also in Table 1 and box $s$.
\end{enumerate}
\end{rem}

\begin{prop}\label{restricciones_nociclicas}
Let $\sigma$, $\tau$, $r$, $s$, $\Bt(\sigma)$, and $\Bt(\tau)$ be as above. Then:
\begin{enumerate}
\item $r$ must be prime and $s=r^j$ with $j\in\Z$.
\item If we call $(\Bt(\sigma)\backslash\Bt(\tau))$ to the subset of all singularities where $\tau$ is locally Cartier and $\sigma$ is not(i.e. singularities where $l_Q\ne 0 = m_Q$), then the  disjoint union $(\Bt(\sigma)\backslash\Bt(\tau))\uplus{\buildrel s \over \dots}\uplus(\Bt(\sigma)\backslash\Bt(\tau))$ %and $(\Bt(\tau)\backslash\Bt(\sigma)\uplus{\buildrel r \over \dots}\uplus(\Bt(\tau)\backslash\Bt(\sigma)$ are 
is in Table 1, box $r$.
\end{enumerate}
\end{prop}

\begin{pf}
The only possibilities that remain (after Table 1 in \cite {C}) are $(r,s)\in\{(2,6),(3,6),(4,4),(4,8),(6,6),(8,8)\}$. We will prove that they are impossible, and this trivially provides 1. 

For $s=6$, we construct the Fano cyclic cover $X_\sigma=\mathrm{Spec}\bigoplus_{i=0}^{r-1}\O_X(i\sigma)$, where we have, as torson divisor, the pullback $\tau'$ of $\tau$. Then there must be a $\Bt$ set of the type of $(\Bt 6,1)$ or $(\Bt6,2)$. First case forces a singularity of type $\frac{1}{12r}$ on $X$ (because there are just one singularity $\frac{1}{12}$ on $X_\sigma$ where $\tau'$ is not Cartier, and so it comes from just one singularity on $X$), which is impossible. If we are dealing with $(\Bt6,2)$, then, by similar reasons, two singularities of type (for $r=3$ it is trivial; for $r=2$, if the two singularities of type $\frac{1}{6}$ in $(\Bt 6.2)$ come from one singularity in $X$, then $\tau$ should have the same shape near both of them), $\frac{1}{6r}$ are forced, and that is also impossible. So $(r,s)=(2,6),(3,6)$ cannot be, and hence $(6,6)$ is also impossible.

For $r=4,8$, it is sufficient to prove that the pair $(4,4)$ is impossible. We think, as in last paragraph, case by case: $(\Bt4.1)$ forces two singularities of type $\frac{1}{32}$ or one of type $\frac{1}{16}$; $(\Bt4.2)$ forces one of type $\frac{1}{48}$; $(\Bt4.3)$ and $(\Bt4.4)$, two of type $\frac{1}{32}$ by the same reason used for $(\Bt6.2)$; finally, $(\Bt4.5)$ would create at least one singularity of type $\frac{1}{16}$ or, else, two of type $\frac{1}{8}$ where $\sigma$ is locally equivalent to $4K_X$ (i.e. $l_Q=4$, which also is not in Table 1.

To prove 2, we know now that $r$ is prime and $s=r^j$ for some $j\in \Z$. Then, when we construct the cyclic cover $X_\tau$ of $X$ related to $\tau$, there will be a torsion divisor $\sigma'$ (the pullback of $\sigma$ by the projection) which will be Cartier on all singularities in the preimage of $\Bt(\tau)$. This is because the order of $\sigma$ near the singularities of $\Bt(\tau)$ is 1 or $r$, and the order of $\tau$ is a nontrivial divisor of $s$, which means that, locally near the singularities of $\Bt(\tau)$, $\sigma$ is a multiple of $\tau$ (whose pullback is Cartier). So $\sigma$ is not Cartier just i the $s$ copies of each singularity in $(\Bt(\sigma)\backslash\Bt(\tau)$. 
\end{pf}

\begin{cor}
In the conditions of Proposition \ref{restricciones_nociclicas}, we have that $r=2$ or $3$, i.e. $(r,s)\in\{ (2,2),(2,4),(2,8),(3,3)\}$. In fact we will see later (Remark \ref{no2,8}) that $(2,8)$ is also impossible.
\end{cor}

\begin{pf}
By Proposition \ref{restricciones_nociclicas}, point 1, $r$ must be prime, and the only prime that we have left from Table 1 in \cite{C} is $5$ (so the group $<\sigma,\tau>$ must be isomorphic to $\Z/(5)\oplus\Z/(5)$). Since no element of Table 1.5 can be considered as the disjoint union of $5$ equal sets, Proposition \ref{restricciones_nociclicas}, point 2 contradicts the possibility of such a subgroup in the Picard of $X$.
\end{pf}

\section{Noncyclic torsion baskets}

Now that we know the form of the torsion subgroup, we can write a table of baskets like Table 1 in \cite{C}. But first of all, we write an example to illustrate the method that we use.

\begin{ex}\label{procedimiento}{\rm
We are searching for a $\Z/(2)\oplus\Z/(4)$ possible torsion basket, so let us suppose that there are two torsion divisors $\sigma$ and $\tau$ of orders 2 and 4 respectively on a Fano threefold $X$. First of all, we use Proposition \ref{restricciones_nociclicas}, point 4 and we observe that the only torsion baskets of order $2$ that can be expressed as a disjoint union of $4$ equal sets are $(\Bt 2.17)$ and $(\Bt2.20)$ that can be found in Table 1.2 in \cite{C}. Then, if we construct the cyclic cover $X_\tau$, we will get a divisor $\sigma'$ of order 2 whose torsion basket is either $(\Bt 2.17)$ or $(\Bt2.20)$. Then, the set of all singularities where $\sigma$ or $\tau$ are not locally Cartier is the union of one singularity of type $\frac{1}{4}(1,1,3)$ (if we had $(\Bt 2.17)$) or two of type $\frac{1}{2}(1,1,1)$ (if we had $(\Bt2.20)$) and $\Bt(\tau)$. Let us consider the possibility of $\Bt(\tau)$ being of type $(\Bt 4.1)$. In the first case, we would deal with a table of values (the entries are the numbers $l_Q$ or $m_Q$ related to the row divisor near the column singularity):
$$\begin{matrix} \ & \vline & \frac{1}{4}(1,1,3) & \frac{1}{4}(1,1,3) & 2\times\frac{1}{8}(1,3,5) \\
\hline \tau & \vline & 0 & 2 & 2 \\
\sigma & \vline & 2 & \mbox{?} & \mbox{?} 
\end{matrix}$$
and we have to choose values for the last two entries of the second row to fit in an element of Table 1.2. But there is not such an element, so this combination is impossible. 

If we try the other option, we would get the table
$$\begin{matrix} & \vline & 2\times\frac{1}{2}(1,1,1) & \frac{1}{4}(1,1,3) & 2\times\frac{1}{8}(1,3,5) \\
\hline \tau & \vline & 0 & 2 & 2 \\
\sigma & \vline & 1 & \mbox{?} & \mbox{?} 
\end{matrix}$$
Now, the only element of Table 1.2 that fits here is $(\Bt 2.14)$ with $\sigma$ as twice the canonical divisor near the singularity of type $\frac{1}{4}(1,1,3)$, cartier near one of the $\frac{1}{8}(1,3,5)$ type singularities and $4K_X$ near the other one. Moreover, there is no contradiction since $\sigma +\tau$ would be associated to $(\Bt4.4)$ and $\sigma+2\tau$ to $(\Bt2.14)$. So this is the only $\Z/(2)\otimes\Z/(4)$ $\Bt$ subset that involve $(\Bt4.1)$.
}\end{ex}

One can now procceed as in Example \ref{procedimiento} and write a table of all possible $\Bt$ subsets of noncyclic torsion groups (any other valid $\Bt$ would be one of these after changing $\tau':=\sigma+\tau$ or something similar):

\vskip 10pt
\centerline{\bf Subsets $\Bt$ for $\Z/(2)\oplus\Z/(2)$:}

\begin{multline*}
(\Bt 2,2.1):={1\over 8}(1,1,7)_{4,0},{1\over 8}(1,1,7)_{4,4},{1\over 8}(1,1,7)_{0,4}
\end{multline*}

\begin{multline*} 
(\Bt 2,2.2):={1\over 8}(1,1,7)_{4,0},{1\over 8}(1,1,7)_{4,4},{1\over 8}(1,3,5)_{0,4} 
\end{multline*} 

\begin{multline*} 
(\Bt 2,2.3):={1\over 8}(1,1,7)_{4,0},{1\over 8}(1,3,5)_{4,4},{1\over 8}(1,3,5)_{0,4} 
\end{multline*} 

\begin{multline*} 
(\Bt 2,2.4):={1\over 8}(1,3,5)_{4,0},{1\over 8}(1,3,5)_{4,4},{1\over 8}(1,3,5)_{0,4} 
\end{multline*} 

\begin{multline*} 
(\Bt 2,2.5):=2\times{1\over 2}(1,1,1)_{1,0},{1\over 4}(1,1,3)_{2,0},{1\over 8}(1,1,7)_{4,4},{1\over 8}(1,1,7)_{0,4} 
\end{multline*} 

\begin{multline*} 
(\Bt 2,2.6):=2\times{1\over 2}(1,1,1)_{1,0},{1\over 4}(1,1,3)_{2,0},{1\over 8}(1,1,7)_{4,4},{1\over 8}(1,3,5)_{0,4} 
\end{multline*} 

\begin{multline*} 
(\Bt 2,2.7):=2\times{1\over 2}(1,1,1)_{1,0},{1\over 4}(1,1,3)_{2,0},{1\over 8}(1,3,5)_{4,4},{1\over 8}(1,3,5)_{0,4} 
\end{multline*} 

\begin{multline*} 
(\Bt 2,2.8):=2\times{1\over 2}(1,1,1)_{1,0},2\times{1\over 2}(1,1,1)_{0,1},{1\over 4}(1,1,3)_{2,0},{1\over 4}(1,1,3)_{0,2},{1\over 8}(1,1,7)_{4,4} 
\end{multline*} 

\begin{multline*} 
(\Bt 2,2.9):=2\times{1\over 2}(1,1,1)_{1,0},2\times{1\over 2}(1,1,1)_{0,1},{1\over 4}(1,1,3)_{2,0},{1\over 4}(1,1,3)_{0,2},{1\over 8}(1,3,5)_{4,4} 
\end{multline*} 

\begin{multline*} 
(\Bt 2,2.10):={1\over 2}(1,1,1)_{1,0}, {1\over 2}(1,1,1)_{0,1}, {1\over 2}(1,1,1)_{1,1}, {1\over 6}(1,1,5)_{3,0}, {1\over 6}(1,1,5)_{0,3}, {1\over 6}(1,1,5)_{3,3} \end{multline*} 

\begin{multline*} 
(\Bt 2,2.11):={1\over 2}(1,1,1)_{0,1}, {1\over 2}(1,1,1)_{1,1}, 2\times{1\over 4}(1,1,3)_{2,0}, {1\over 6}(1,1,5)_{0,3}, {1\over 6}(1,1,5)_{3,3} \end{multline*} 

\begin{multline*} 
(\Bt 2,2.12):=4\times{1\over 2}(1,1,1)_{1,0}, {1\over 2}(1,1,1)_{0,1}, {1\over 2}(1,1,1)_{1,1}, {1\over 6}(1,1,5)_{0,3}, {1\over 6}(1,1,5)_{3,3} 
\end{multline*} 

\begin{multline*} 
(\Bt 2,2.13):={1\over 2}(1,1,1)_{1,1}, 2\times{1\over 4}(1,1,3)_{2,0}, 2\times{1\over 4}(1,1,3)_{0,2}, {1\over 6}(1,1,5)_{3,3} 
\end{multline*}

\begin{multline*} 
(\Bt 2,2.14):=4\times{1\over 2}(1,1,1)_{1,0}, {1\over 2}(1,1,1)_{0,1}, 2\times{1\over 4}(1,1,3)_{2,2}, {1\over 6}(1,1,5)_{0,3} 
\end{multline*} 

\begin{multline*} 
(\Bt 2,2.15):=2\times{1\over 4}(1,1,3)_{2,0}, 2\times{1\over 4}(1,1,3)_{0,2}, 2\times{1\over 4}(1,1,3)_{2,2} 
\end{multline*} 

\begin{multline*} 
(\Bt 2,2.16):=4\times{1\over 2}(1,1,1)_{1,0}, 2\times{1\over 4}(1,1,3)_{0,2}, 2\times{1\over 4}(1,1,3)_{2,2} 
\end{multline*} 

\begin{multline*} 
(\Bt 2,2.17):=4\times{1\over 2}(1,1,1)_{1,0}, 4\times{1\over 2}(1,1,1)_{0,1}, {1\over 2}(1,1,1)_{1,1}, {1\over 6}(1,1,5)_{3,3} 
\end{multline*} 

\begin{multline*} 
(\Bt 2,2.18):=2\times{1\over 2}(1,1,1)_{1,0}, 2\times{1\over 2}(1,1,1)_{0,1}, \\ 2\times{1\over 2}(1,1,1)_{1,1}, {1\over 4}(1,1,3)_{2,0}, {1\over 4}(1,1,3)_{0,2}, {1\over 4}(1,1,3)_{2,2} 
\end{multline*} 

\begin{multline*} 
(\Bt 2,2.19):=4\times{1\over 2}(1,1,1)_{1,0}, 4\times{1\over 2}(1,1,1)_{0,1}, 2\times{1\over 4}(1,1,3)_{2,2} 
\end{multline*} 

\begin{multline*} 
(\Bt 2,2.20):=4\times{1\over 2}(1,1,1)_{1,0}, 4\times{1\over 2}(1,1,1)_{0,1}, 4\times{1\over 2}(1,1,1)_{1,1} 
\end{multline*}

\vskip.2in
\centerline{\bf Subsets $\Bt$ for $\Z/(2)\oplus\Z/(4)$:}

\begin{multline*} 
(\Bt 2,4.1):=2\times{1\over 2}(1,1,1)_{1,0},{1\over 4}(1,1,3)_{2,2},{1\over 8}(1,3,5)_{4,2},{1\over 8}(1,3,5)_{0,2} 
\end{multline*} 

\begin{multline*} 
(\Bt 2,4.2):=2\times{1\over 2}(1,1,1)_{1,0}, 2\times{1\over 2}(1,1,1)_{1,1}, {1\over 4}(1,1,3)_{2,1}, \\ {1\over 4}(1,1,3)_{0,1}, {1\over 4}(1,1,3)_{2,3}, {1\over 4}(1,1,3)_{0,3} 
\end{multline*}

\vskip.2in
\centerline{\bf Subsets $\Bt$ for $\Z/(3)\oplus\Z/(3)$:}

\begin{multline*} 
(\Bt 3,3.1):={1\over 3}(1,1,2)_{1,0}, {1\over 3}(1,1,2)_{2,0}, {1\over 3}(1,1,2)_{1,1}, \\ {1\over 3}(1,1,2)_{2,1}, {1\over 3}(1,1,2)_{1,2}, {1\over 3}(1,1,2)_{2,2}, \\ {1\over 3}(1,1,2)_{0,1}, {1\over 3}(1,1,2)_{0,2} 
\end{multline*}

\begin{rem}\label{no2,8}{\rm
It is easy to use the method developed in Example \ref{procedimiento} to see that there does not exist a possible $\Z/(2)\oplus\Z/(8)$ subgroup.
}\end{rem}

\begin{rem}{\rm
After this, we can repeat what we have done to check if we can find three torsion divisors $\sigma,\tau,\upsilon$ of orders $r,s$ and $t$ where $r|s$ and $s|t$. We can repeat the reasoning to see that:
\begin{itemize}
\item $r$ must be prime with $s=r$, $t=r^j$, $j\in \Z$.
\item extending naturally the above notation, we get $(\Bt(\sigma)\backslash(\Bt(\tau)\cup\Bt(\upsilon)))\uplus{\buildrel st \over \dots}\uplus(\Bt(\sigma)\backslash(\Bt(\tau)\cup\Bt(\upsilon)))$ must be a valid $\Bt$ subset of order $r$
\end{itemize}

This implies that $r=s=2$ (because no $\Bt$ subset of order $3$ can be expressed as a disjoint union of $9$ equal sets) and we can now write the corresponding table for three torsion divisors using the natural generalization of the method in Example \ref{procedimiento}.
}\end{rem}

\vskip 10pt
\centerline{\bf Subsets $\Bt$ for $\Z/(2)\oplus\Z/(2)\oplus\Z/(2)$:}

\begin{multline*} 
(\Bt 2,2,2.1):={1\over 4}(1,1,3)_{2,0,0}, {1\over 4}(1,1,3)_{2,2,0}, {1\over 4}(1,1,3)_{2,2,2}, \\ {1\over 4}(1,1,3)_{2,0,2}, {1\over 4}(1,1,3)_{0,2,0}, {1\over 4}(1,1,3)_{0,2,2}, {1\over 4}(1,1,3)_{0,0,2}  
\end{multline*} 

\begin{multline*} 
(\Bt 2,2,2.2):=2\times{1\over 2}(1,1,1)_{1,0,0}, 2\times{1\over 2}(1,1,1)_{1,1,0}, 2\times{1\over 2}(1,1,1)_{0,1,0}, \\ {1\over 4}(1,1,3)_{2,0,2}, {1\over 4}(1,1,3)_{2,2,2}, {1\over 4}(1,1,3)_{0,2,2}, {1\over 4}(1,1,3)_{0,0,2}  
\end{multline*} 

\begin{multline*} 
(\Bt 2,2,2.3):=2\times{1\over 2}(1,1,1)_{1,0,0}, 2\times{1\over 2}(1,1,1)_{1,0,1}, 2\times{1\over 2}(1,1,1)_{0,1,0}, \\ 2\times{1\over 2}(1,1,2)_{0,1,1}, {1\over 4}(1,1,3)_{2,2,0}, {1\over 4}(1,1,3)_{2,2,2}, {1\over 4}(1,1,3)_{0,0,2}  
\end{multline*} 

\begin{multline*} 
(\Bt 2,2,2.4):=2\times{1\over 2}(1,1,1)_{1,0,0}, 2\times{1\over 2}(1,1,1)_{1,0,1}, 2\times{1\over 2}(1,1,1)_{1,1,0}, \\ 2\times{1\over 2}(1,1,2)_{1,1,1}, 2\times{1\over 2}(1,1,2)_{0,1,0}, 2\times{1\over 2}(1,1,2)_{0,1,1}, 2\times{1\over 2}(1,1,2)_{0,0,1}  
\end{multline*}

\begin{rem}{\rm
We check now what happens for the case of four generators and we can only consider $\Z/(2)\oplus\Z/(2)\oplus\Z/(2)\oplus\Z/(2)$ since no $\Bt$ of order $r\ne 2$ can be expressed as a disjoint union of $r^ir^jr^k$ equal sets with $i,j,k>0$. We can this time get just one basket.
}\end{rem}

\vskip 10pt
\centerline{\bf Subsets $\Bt$ for $\Z/(2)\oplus\Z/(2)\oplus\Z/(2)\oplus\Z/(2)$:}

\begin{multline*} 
(\Bt 2,2,2.4):={1\over 2}(1,1,1)_{1,0,0,0}, {1\over 2}(1,1,1)_{1,1,0,0}, {1\over 2}(1,1,1)_{1,1,0,1}, \\ {1\over 2}(1,1,1)_{1,1,1,0}, {1\over 2}(1,1,2)_{1,1,1,1}, {1\over 2}(1,1,2)_{1,0,1,0}, {1\over 2}(1,1,2)_{1,0,1,1}, \\ {1\over 2}(1,1,2)_{1,0,0,1}, {1\over 2}(1,1,1)_{0,1,0,0}, {1\over 2}(1,1,1)_{0,1,0,1}, {1\over 2}(1,1,1)_{0,1,1,0}, \\ {1\over 2}(1,1,2)_{0,1,1,1}, {1\over 2}(1,1,2)_{0,0,1,0}, {1\over 2}(1,1,2)_{0,0,1,1}, {1\over 2}(1,1,2)_{0,0,0,1} 
\end{multline*}

\begin{rem}
We can now think as before and see that, since no $\Bt$ of order $2$ can be expressed as the disjoint union of 16 equal sets, it is impossible to find five independent torsion divisors.
\end{rem}

\section{Finding quotients}

Now we know all possible shapes of a distinct subset of the basket of a noncyclic Fano-Enriques. We can now write all Fano-Enriques threefolds with noncyclic torsion that can be obtained as quotients of the Fano threefolds in Reid's, Fletcher's and Alt{\i}nok's lists. The way we work is the same as in \cite{C}, with the obvious modifications:
\begin{itemize}
\item We can calculate the basket of a noncyclic covering just repeating the method for a cyclic one as many times as needed.
\item We now use the graded ring:
$$\bigoplus_{n\in\Z, a_i\in \Z/(r_i)}H^0(-nK_X+a_1\sigma_1+...+a_l\sigma_l)$$
where $\sigma_1,...,\sigma_l$ are the generators of the torsion subgroup of Pic$X$ and $r_1,...,r_l$ are their orders.
\item To check that all singularities are as we wish, we just have to check them cyclic action by cyclic action (in the proper order). 
\end{itemize}
It will probably be clearer in the next example, which is analog to the (split in parts through the paper) example in \cite{C} of a $\Z/(5)$ quotient:

\begin{ex}{\rm
Let us consider the intersection $Y$ of three quadrics in $\P^6$. If we want to find a $G:=\Z/(2)\oplus\Z/(4)$ quotient, we can only take as torsion subset ($\Bt 2,4.2$), since the two singularitiesof type $\frac{1}{8}$ in ($\Bt 2,4.1$) come from two singularities of type $\frac{1}{2}$ in the cover (that is because none of the torsion divisors generate the class divisor group of a sufficiently small neighbourhood of the singularity). If $X$ is a quotient of $Y$ by a $G$ action, then the selfintersection of the anticanonical bundle must be $\frac{8}{ord(G)}=1$. We have now all the data we need to compute the Hilbert series of the graded ring of $X$ (in this case, this ring is graded by $\Z\oplus G$):
$$\bigoplus_{n\in\Z,a\in\Z/(2),b\in \Z/(4)} H^0(X, \O(-nK_X+a\sigma+b\tau)),$$
which is
$$HS:=\sum_{n\in\Z,a\in\Z/(2),b\in \Z/(4)} h^0(X, \O(-nK_X+a\sigma+b\tau))t^ne_1^a,e_2^b.$$
Where, of course, $a^2=b^4=1$.

We use Alt\i nok's formula for the plurigenus (see \cite{ABR} or \cite{C}) and torsion formula (\cite[Lemma 3.2]{C}) to get the series $\sum_n h^0(X, \O(-nK_X+a\sigma+b\tau))t^ne_1^a,e_2^b$ for fixed $a$ and $b$. The result is:
\begin{multline*}
1 + t + 3t^2 + 7t^3 + 17t^4 + 29t^5 + 47t^6 + 71t^7 + 105t^8+...\\
+e_2( t + 3t^2 + 8t^3 + 16t^4 + 29t^5 + 47t^6 + 72t^7 + 104t^8 +...)+\\
+e_2^2(4t^2 + 8t^3 + 16t^4 + 28t^5 + 48t^6 + 72t^7 + 104t^8 + 144t^9 + ...)+\\
+e_2^3( t + 3t^2 + 8t^3 + 16t^4 + 29t^5 + 47t^6 + 72t^7 + 104t^8 + ...)+\\
e_1( t + 3t^2 + 8t^3 + 16t^4 + 29t^5 + 47t^6 + 72t^7 + 104t^8 + ...)+\\
e_1e_2( t + 3t^2 + 8t^3 + 16t^4 + 29t^5 + 47t^6 + 72t^7 + 104t^8 + ...)+\\
e_1e_2^2( t + 3t^2 + 8t^3 + 16t^4 + 29t^5 + 47t^6 + 72t^7 + 104t^8 + ...)+\\
e_1e_2^3( t + 3t^2 + 8t^3 + 16t^4 + 29t^5 + 47t^6 + 72t^7 + 104t^8 + ...)
\end{multline*}
Therefore, we need generators in tridegrees $(1,+,+),(1,+,i),(1,+,-i),(1,-,+),(1,-,i),(1,-,-),(1,-,-i)$ (where $i^2=-1$). In fact, we have the hilbert series is:
$$\frac{1 - 2t^2 + t^4 + e_2^2(-t^2 + 2t^4 - t^6) }{(1-t)(1-e_2t)(1-e_2^3t)(1-e_1t)(1-e_1e_2t)(1-e_1e_2^2t)(1-e_1e_2^3t)}$$
This means that $\Z/(2)\oplus\Z/(4)$ acts in $\P^6$ by $(x_0:x_1:x_2:x_3:x_4:x_5:x_6)\mapsto(x_0:x_1:x_2:-x_3:-x_4:-x_5:-x_6)$ (the first generator) and $(x_0:x_1:x_2:x_3:x_4:x_5:x_6)\mapsto(x_0:ix_1:-ix_2:x_3:ix_4:-x_5:-ix_6)$ (the secon one). On the other side, we know by the numerator that there are two invariant equations (so they are of type $c_1x_0^2+c_2x_1x_2+c_3x_3+c_4x_4x_6+c_5x_5^2$) and the other one is invariant for the first map and multiplied by $-1$ by the second one (so of shape $c_2x_1^2+c_2x_2^2+c_3x_3x_5+c_4x_4^2+c_5x_6^2$). It is easy to check that the general element is quasi smooth and the non-\'etal\'e points are as desired.
}\end{ex}

\section{Lists of noncyclic Fano--Enriques}\label{listas}
\vskip 10pt
\centerline{\bf Noncyclic Fano--Enriques threefolds from codimension 2 Fano threefolds:}
\vskip 10pt
\leftline{\bf No. 1}
\begin{itemize}
\item {\bf cover:} $Y_{4,4}\subset\P(1,1,1,2,2,2)$
\item {\bf action:} $\Z/(2)\oplus\Z/(2)$ acts by $([+,-],[-,+],[-,-],[+,-],\ [-,+],[-,-])$
\item $\Bt=(\Bt 2,2.20)$ 
\item $\B \backslash \Bt={1\over 2}(1,1,1)$
\item Both equations are invariant by the action.
\end{itemize}

\vskip 10pt
\centerline{\bf Noncyclic Fano--Enriques threefolds from codimension 3 Fano threefolds:}
\vskip 10pt
\leftline{\bf No. 1a}
\begin{itemize}
\item {\bf cover:} $Y_{2,2,2}\subset\P(1,1,1,1,1,1,1)$
\item {\bf action:} $\Z/(2)\oplus\Z/(2)$ acts by $([+,+],[+,-],[+,-],[-,+],[-,+],[-,-],[-,-])$
\item $\Bt=(\Bt 2,2.20)$ 
\item $\B \backslash \Bt=\emptyset$
\item All equations are invariant by the action.
\end{itemize}
\vskip 10pt
\leftline{\bf No. 1b}
\begin{itemize}
\item {\bf cover:} $Y_{2,2,2}\subset\P(1,1,1,1,1,1,1)$
\item {\bf action:} $\Z/(2)\oplus\Z/(4)$ acts by $([+,+],[+,i],[+,-i],[-,+],[-,i],[-,-],[-,-i])$
\item $\Bt=(\Bt 2,4.2)$ 
\item $\B \backslash \Bt=\emptyset$
\item Two equations are invariant by the action and the other one has second degree $[0,2]$.
\end{itemize}
\vskip 10pt
\leftline{\bf No. 1c}
\begin{itemize}
\item {\bf cover:} $Y_{2,2,2}\subset\P(1,1,1,1,1,1,1)$
\item {\bf action:} $\Z/(2)\oplus\Z/(2)\oplus\Z/(2)$ acts by $([+,+,-],[+,-,+],[+,-,+],[-,+,+],[-,+,-],[-,-,+],[-,-,-])$
\item $\Bt=(\Bt 2,2,2.4)$ 
\item $\B \backslash \Bt=\emptyset$
\item All equations are invariant by the action.
\end{itemize}

\end{document}